\documentclass{amsart}%
\usepackage{amsfonts}
\usepackage{amsmath}
\usepackage{amssymb}
\usepackage{graphicx}%
\providecommand{\U}[1]{\protect\rule{.1in}{.1in}}
\newtheorem{theorem}{Theorem}
\theoremstyle{plain}

\newtheorem{definition}{Definition}
\newtheorem{example}{Example}

\newtheorem{lemma}{Lemma}

\newtheorem{problem}{Problem}
\newtheorem{proposition}{Proposition}
\newtheorem{remark}{Remark}
\newtheorem{solution}{Solution}

\numberwithin{equation}{section}
\begin{document}
\title{Reconstruction of the one-dimensional thick distribution theory }
\author{Yunyun Yang}
\address{Department of Mathematics, Hefei University of Technology}
\email{yangyunyun@hfut.edu.cn}
\urladdr{http://maths.hfut.edu.cn/2016/1126/c7700a213853/page.htm}
\thanks{Thanks to Hefei University of Technology to support this research. The grant
number from Hefei University of Technology is 407-0371000086}
\date{August 20, 2019}
\keywords{thick test functions, thick distributions, thick delta functions, Hadmard
finite part, asymptotic expansion}

\begin{abstract}
The theory of thick distributions (both in dimension 1 and in higher
dimensions) was constructed in recent years \cite{Thick1, Thick}. However this
theory of distributions with one thick point in dimension one is very
different from that in higher dimensions. In this paper the author uses the
language of asymptotic analysis to reconstruct the 1-dimensional thick
distribution theory and to incorporate it into the framework of the
higher-dimensional thick distribution theory. Some new concepts and
interesting results appear in this paper from viewing singular functions in a
different way.

\end{abstract}
\maketitle





\section{Introduction}

In recent years, some research in mathematics and physics suggests us to
consider non-smooth functions as test functions in the theory of
distributions. In particular, one considers functions with a point
singularity. The emergence of the special point either corresponds to the
singularity of the field equations, or to the non-linearity caused by the
multiplication of distributions \cite{Blanchet, BlaFay, Blinder, Bowen,
Franklin, Gsponer, Thick}. Classical distribution theory cannot handle these problems.

In the past decades, researchers have done some work to solve these
problems. Blanchet and Faye develop such a scheme in the context of finite
parts, pseudo-functions and Hadamard regularization, to study the dynamics of
point particles in high post-Newtonian approximations of general relativity
\cite{BlaFay}. Their work was built upon the the discussion by Sellier of 
Hadamard finite part of smooth functions with one point singularity
\cite{Sellier1}.

Considering distributions with a point singularity can also be
applied to the research of point-source fields. Basic electrodynamics tells us
that, because of $\nabla^{2}\frac{1}{r}=-4\pi\delta\left(  \mathbf{x}\right)
$ and the Poisson's equation $\nabla^{2}\Phi\left(  \mathbf{x}\right)
=-\rho\left(  \mathbf{x}\right)  /\epsilon_{0},$ we can describe the potential
of a point-source field caused by an electron $\rho\left(  \mathbf{x}\right)
=q\delta\left(  \mathbf{x}\right)  .$ Yet applying the Laplace operator
$\nabla^{2}=\frac{1}{r^{2}}\frac{\partial}{\partial r}r^{2}\frac{\partial
}{\partial r}+...$directly on $\frac{1}{r}$ cannot give us the desired
$\nabla^{2}\frac{1}{r}=-4\pi\delta\left(  \mathbf{x}\right)  .$ In order to
solve this, Blinder \cite{Blinder} and Hu \cite{Hu} suggest to use $sgn\left(
r\right)  /r$ instead of $1/r$ to have $\nabla^{2}\frac{sgn\left(  r\right)
}{r}=-2\delta\left(  r\right)  /r^{2}=-4\pi\delta\left(  \mathbf{x}\right)  .$
Yet $sgn\left(  r\right)  /r$ and $-2\delta\left(  r\right)  /r^{2}$ both are
NOT\ well-defined distributions. In fact, they are so called ``thick distribution''.

Paskusz \cite{Paskusz} points out, it is not rigorous to say $H\left(
x\right)  \delta\left(  x\right)  =\frac{1}{2}\delta\left(  x\right)  $ based
on $H\left(  x\right)  =H^{2}\left(  x\right)  .$ Indeed, taking
distributional derivative on both sides of $H\left(  x\right)  =H^{2}\left(
x\right)  $ does give us $2H\left(  x\right)  \delta\left(  x\right)
=\delta\left(  x\right)  .$ Yet, if one continues to multiply both sides with
$2H\left(  x\right)  $ one gets $\delta\left(  x\right)  =2H\left(  x\right)
\delta\left(  x\right)  =4H^{2}\left(  x\right)  \delta\left(  x\right)
=4H\left(  x\right)  \delta\left(  x\right)  ,$ and then gets $2=4$. In fact,
both $H^{2}\left(  x\right)  $ and $H\left(  x\right)  \delta\left(  x\right)
$ are not well-defined distributions. One cannot simply apply distributional
derivatives on both sides of the above equation.

Bowen, \cite{Bowen}, points out that, when we multiply two distributions, the
product rule may not apply. For example, let $n_{j_k}=x_{j_k}/r,$ where $x_{j_k}$ is the $j_k$-th coordinate and $r$ is the radius $r=\sqrt{x_1^2+...+x_n^2}$. We have
\[
\frac{\overline{\partial}}{\partial x_{i}}\left(  \frac{n_{j_{1}}n_{j_{2}%
}n_{j_{3}}}{r^{2}}\right)  \neq n_{j_{1}}n_{j_{2}}\frac{\overline{\partial}%
}{\partial x_{i}}\left(  \frac{n_{j_{3}}}{r^{2}}\right)  +\frac{n_{j_{3}}%
}{r^{2}}\frac{\overline{\partial}}{\partial x_{i}}\left(  n_{j_{1}}n_{j_{2}%
}\right)  .
\]
where the $\overline{\partial}$ means the distributional derivative. However,
the equation holds when one considers "thick distributional derivative".

A one-dimensional thick distribution theory was first developed by Estrada and
Fulling in order to solve some of these dilemmas \cite{Thick1}. They define
the "thick test function space" as the topological vector space consists of
all compactly supported functions which are smooth on $\mathbb{R}%
\backslash\left\{  a\right\}  $, and whose one-sided derivatives at $x=a$
exist. They defined the "thick distribution space" as the dual space of the
thick test function space.

The thick distribution theory of higher dimensions was later developed by
Estrada and Yang \cite{Thick}. We also published a series of papers to further
develop the application of thick distributions \cite{Yang2,Yang4, Yang1, Yang3}. In higher dimensions, we define the "thick test functions" as those
compactly supported functions that are smooth on $\mathbb{R}^{n}%
\backslash\left\{  \mathbf{a}\right\}  ,$ and which have a strong asymptotic
expansion near $\mathbf{a}:\phi\left(  \mathbf{a}+\mathbf{x}\right)  \sim
\sum_{j=m}^{\infty}a_{j}\left(  \mathbf{w}\right)  r^{j}$ when $\mathbf{x}%
\rightarrow\mathbf{0,}$ where $a_{j}\left(  \mathbf{w}\right)  $ is a smooth
function on the unit sphere $\mathbb{S}^{n-1}.$ Given a proper topology, the
space of all such functions is a topological vector space. We then define the
dual space as the "thick distribution space".

One could see that because $\mathbb{R}^{n}$ ($n\geq2)$ is connected when
taking away a single point $\mathbf{a}$, while $\mathbb{R}$ is disconnect
doing the same thing, the construction of the higher dimensional thick
distribution theory is very different from that of the one-dimensional case.
One could not talk about "one-sided derivatives" in the higher dimensional
case. Instead, we used the asymptotic expansion behavior near the singularity
to describe thick test functions.

The purpose of this present article is to try to combine these two seemingly
different theories, and to incorporate the one-dimensional thick distribution
theory into the framework of the higher dimensional thick distribution theory.
This work will give us some insights of the function space itself.

About the notation: $d/dx$ means the ordinary derivative, $\overline{d}/dx$
means the distributional derivative, $d^{\ast}/dx$ means the thick
distributional derivative.

This is the structure of this article: in Section \ref{sec2} we give a brief
review of the established thick distribution theory. Section \ref{sec3} is the
main part of this article, we introduce a new one-dimensional thick
distribution theory there, which is more allied with the higher dimensional
thick distribution theory, and the established one-dimensional thick
distribution theory is a special case of the new thick distribution
theory. In this section, we firstly identify $\mathbb{R}\backslash\left\{
0\right\}  $ with $\mathbb{S}^{0}\times\mathbb{R}_{+},$ where "$\mathbb{S}
^{0}"$ denotes the "0-dimensional unit sphere", namely, the 2 end-points of
the line segment $[-1,1].$ Under this viewpoint, we define the meaning of a
function $f\left(  x\right)  =f\left(  \mathbf{w},r\right)  $ having an
asymptotic expansion $\sum_{j=m}^{\infty}a_{j}\left(  \mathbf{w}\right)
r^{j}$ as $r\rightarrow0^{+},$ where $\mathbf{w}\in\mathbb{S}^{0}.$ Then, we 
give a few examples where viewing $f\left(  x\right)  $ as $f\left(
\mathbf{w},r\right)  $ are quite elegant, as well as a few lemmas discussing a
few properties of $f\left(  \mathbf{w},r\right)  .$ In subsection
\ref{subsec3.2} we introduce a new thick test function space in
one-dimensional case. Under our definition, the established thick test
function space in one-dimensional case is a closed subspace. In subsections
\ref{subsec3.3} and \ref{subsec3.4}, we reconstruct the one-dimensional thick
distribution theory, and provide a few examples of the new thick
distributions. At the end of the article we provide a solution of a problem
pointed out by Paskusz, of multiplying Schwartz distributions, in the context
of the new one-dimensional thick distributoin theory.

\section{\bigskip Review of the established 1-dim and higher dimensional thick
distribution theory\label{sec2}}

\subsection{One-dimensional thick distribution theory}

One-dimensional thick distribution theory was developed by Estrada and Fulling
in 2007 \cite{Thick1}. In their paper, the "thick test function space" is
defined as the topological vector space consists of all compactly supported
functions which are smooth on $\mathbb{R}\backslash\left\{  a\right\}  $, and
whose one-sided derivatives at $x=a$ exist. One denotes thick test function
space $\mathcal{D}_{\ast,a}\left(  \mathbb{R}\right)  .$ Given a proper
topology, the usual test function space $\mathcal{D}\left(  \mathbb{R}\right)
$ is a closed subspace of $\mathcal{D}_{\ast,a}\left(  \mathbb{R}\right)  .$
Then they define the dual space of $\mathcal{D}_{\ast,a}\left(  \mathbb{R}%
\right)  $ as ``thick distribution space", denoted $\mathcal{D}_{\ast
,a}^{\prime}\left(  \mathbb{R}\right)  .$ If $a=0,$ we simply denote the space
of thick test functions and thick distributions $\mathcal{D}_{\ast}\left(
\mathbb{R}\right)  $ and $\mathcal{D}_{\ast}^{\prime}\left(  \mathbb{R}%
\right)  $, respectively.

By the Hahn-Banach theorem, the usual distribution space $\mathcal{D}^{\prime
}\left(  \mathbb{R}\right)  $ is just a projection of $\mathcal{D}_{\ast
,a}^{\prime}\left(  \mathbb{R}\right)  .$ Namely, any usual distribution $f$
has a lifting $\widetilde{f}$ in $\mathcal{D}_{\ast,a}^{\prime}\left(
\mathbb{R}\right)  ,$ such that $\pi\left(  \widetilde{f}\right)  =f,$ where
$\pi$ is the projection operator.

An example of thick distribution in $\mathcal{D}_{\ast}^{\prime}\left(
\mathbb{R}\right)  $ would be $\delta_{+}\left(  x\right)  :$ for any thick test function $\phi\left(  x\right)  $, we have
\[
\left\langle \delta_{+}\left(  x\right)  ,\phi\left(  x\right)  \right\rangle
=\phi_{+}\left(  0\right)  ,
\]
where $\phi_{+}\left(  0\right)  $ denotes $\underset{x\rightarrow0^{+}}{\lim
}\phi\left(  x\right)  ,$ the rightside limit of $\phi\left(  x\right)  ;$ and
$\phi_{-}\left(  0\right)  $ denotes $\underset{x\rightarrow0^{-}}{\lim}%
\phi\left(  x\right)  .$ Moreover, one can define an extension of the Dirac
delta function as:
\begin{equation}
\delta_{\ast,\lambda}\left(  x\right)  =\lambda\delta_{+}\left(  x\right)
+\left(  1-\lambda\right)  \delta_{-}\left(  x\right)
.\label{thick delta old}%
\end{equation}

Namely, the projection of $\delta_{\ast,\lambda}\left(  x\right)  $ onto the
usual distribution space is the Dirac delta function $\delta\left(  x\right)
$: if we use $i:\mathcal{D}\left(  \mathbb{R}\right)  \hookrightarrow
\mathcal{D}_{\ast}\left(  \mathbb{R}\right)  $ to denote the inclusion and
$\pi:\mathcal{D}_{\ast}^{\prime}\left(  \mathbb{R}\right)  \rightarrow
\mathcal{D}^{\prime}\left(  \mathbb{R}\right)  $ the projection, then one
clearly sees that, for any $\phi\left(  x\right)  \in\mathcal{D}\left(
\mathbb{R}\right)  ,$ we have $\phi_{+}\left(  0\right)  =\phi_{-}\left(
0\right)  =\phi\left(  0\right)  ,$ then,
\begin{align*}
\left\langle \pi\left(  \delta_{\ast,\lambda}\left(  x\right)  \right)
,\phi\left(  x\right)  \right\rangle  &  =\left\langle \delta_{\ast,\lambda
}\left(  x\right)  , i\left(  \phi\left(  x\right)  \right)  \right\rangle \\
&  =\lambda\phi_{+}\left(  0\right)  +\left(  1-\lambda\right)  \phi
_{-}\left(  0\right)  \\
&  =\phi\left(  0\right) =\left\langle \delta\left(  x\right)  ,\phi\left(
x\right)  \right\rangle .
\end{align*}

In particular, if $\lambda=\frac{1}{2},$%
\[
\delta_{\ast,1/2}\left(  x\right)  =\frac{1}{2}\delta_{+}\left(  x\right)
+\frac{1}{2}\delta_{-}\left(  x\right)  .
\]

Let $\mathcal{S}_{\ast}\left(  \mathbb{R}\right)  $ denote the space of smooth
functions with a thick point at the origin, and with fast decay at infinity.
Then in \cite{Thick1} Estrada and Fulling proved that the space of Fourier
transform of $\mathcal{S}_{\ast}\left(  \mathbb{R}\right)  $ would be the
space $\mathcal{W}\left(  \mathbb{R}\right)  .$ Namely, the space consists of
those functions $\psi\in C^{\infty}\left(  \mathbb{R}\right)  $ that admits an
asymptotic expansion of the type%
\[
\psi\left(  x\right)  \sim\underset{n=1}{\overset{\infty}{\sum}}c_{n}%
x^{-n},\text{ \ \ \ as }\left\vert x\right\vert \rightarrow\infty,
\]
for some constants $c_{1},c_{2},c_{3},...$. The correponding dual space
$\mathcal{W}^{\prime}$ is therefore the Fourier transform of the tempered
thick distributions $\mathcal{S}_{\ast}^{\prime}.$

\subsection{Higher dimensional thick distribution theory}

In paper \cite{Thick}, Yang and Estrada developed a theory of thick
distributions in dimension $n\geq2.$ There is a essential difference between
the topology of $\mathbb{R}$ and that of $\mathbb{R}^{n}:$ taking away a point
in $\mathbb{R}$ will make it into two disconnected parts while doing the same
thing for $\mathbb{R}^{n}$ will not. Hence we cannot use the so called "jump
discontinuity" to describle singularities of functions in $\mathbb{R}^{n}.$
Therefore, the theory of thick distributions in higher dimensions is quite
different from the established one-dimensional theory.

\begin{definition}
Let $\phi$ be defined in $\mathbb{R}^{n}\backslash\left\{  \mathbf{0}\right\}
.$ We say that $\phi$ has the asymptotic expansion $\sum_{j=m}^{\infty}%
a_{j}\left(  \mathbf{w}\right)  r^{j}$ as $\mathbf{x}\rightarrow\mathbf{0}$ if
for all $M\geq m,M\in\mathbb{Z}$,
\begin{equation}
\underset{r\rightarrow0^{+}}{\lim}\left\vert \phi\left(  \mathbf{x}\right)
-\underset{j=m}{\overset{M}{\sum}}a_{j}\left(  \mathbf{w}\right)  r^{j}\right\vert
r^{-M}=0,\text{ \ \ uniformly on }\mathbf{w}\in\mathbb{S}\text{.}
\label{asy expansion}%
\end{equation}
In this case we write $\phi(\mathbf{x})\sim\sum_{j=m}^{\infty}a_{j}\left(  \mathbf{w}%
\right)  r^{j}$ as $\mathbf{x}\rightarrow\mathbf{0}.$
\end{definition}

Let $\mathbf{p}\in\mathbb{N}^{n}$ denote a multi-index and $\left(
\partial/\partial\mathbf{x}\right)  ^{\mathbf{p}}=\left(  \partial^{\left\vert
\mathbf{p}\right\vert }\right)  /\partial x_{1}^{p_{1}}...\partial
x_{n}^{p_{n}},\left\vert \mathbf{p}\right\vert =p_{1}+...+p_{n}.$ If
$a_{j}\left(  \mathbf{w}\right)  \in C^{\infty}\left(  \mathbb{S}%
^{n-1}\right)  $ is a smooth function on the unit sphere $\mathbb{S}^{n-1},$ then applying $\left(  \partial/\partial\mathbf{x}\right)  ^{\mathbf{p}}$ on
$a_{j}\left(  \mathbf{w}\right)  r^{j}$ we will obtain another homogeneous
function of degree $j-\left\vert \mathbf{p}\right\vert ,$ which we denote
$a_{j-\left\vert \mathbf{p}\right\vert ,\mathbf{p}}\left(  \mathbf{w}\right)
r^{j-\left\vert \mathbf{p}\right\vert }.$ Hence formally applying $\left(
\partial/\partial\mathbf{x}\right)  ^{\mathbf{p}}$ on $\sum_{j=m}^{\infty
}a_{j}\left(  \mathbf{w}\right)  r^{j}$ will give us another expansion of the
form $\sum_{j=m-\left\vert \mathbf{p}\right\vert }^{\infty}a_{j,\mathbf{p}%
}\left(  \mathbf{w}\right)  r^{j}.$ In general, asymptotic expansions cannot be
differentiated \cite{Greenbook}; if the asymptotic expansion of the
differentiation of the function is the same as the term-by-term
differentiation of the expansion, then the expansion is called
"\textsl{strong}".

\begin{definition}
Let $\phi\in C^{\infty}\left(  \mathbb{R}^{n}\backslash\left\{  \mathbf{0}%
\right\}  \right)  .$We say that the expansion $\phi\left(  \mathbf{x}\right)
\sim\sum_{j=m}^{\infty}a_{j}\left(  \mathbf{w}\right)  r^{j}$ as
$\mathbf{x}\rightarrow0$ is strong if for each $\mathbf{p}\in\mathbb{N}^{n}$
the asymptotic development of $\left(  \partial/\partial\mathbf{x}\right)
^{\mathbf{p}}\phi\left(  \mathbf{x}\right)  $ as $\mathbf{x}\rightarrow
\mathbf{0}$ exists and equals $\sum_{j=m-\left\vert \mathbf{p}\right\vert
}^{\infty}a_{j.\mathbf{p}}\left(  \mathbf{w}\right)  r^{j},$ the term-by-term
differentiation of $\sum_{j=m}^{\infty}a_{j}\left(  \mathbf{w}\right)  r^{j}.$
\end{definition}

\begin{definition}
Let $\mathcal{D}_{\ast,\mathbf{a}}\left(  \mathbb{R}^{n}\right)  $ denote the
vector space of all smooth functions $\phi$ defined in $\mathbb{R}%
^{n}\backslash\left\{  \mathbf{a}\right\}  ,$ with support of the form
$K\backslash\left\{  \mathbf{a}\right\}  ,$ where $K$ is compact in
$\mathbb{R}^{n},$ that admits a strong asymptotic expansion of the form
\begin{equation}
\phi\left(  \mathbf{a}+\mathbf{x}\right)  =\phi\left(  \mathbf{a}%
+r\mathbf{w}\right)  \sim\sum_{j=m}^{\infty}a_{j}\left(  \mathbf{w}\right)
r^{j},\text{ \ \ as }\mathbf{x}\rightarrow\mathbf{0},\text{ }
\label{test function higher dim}%
\end{equation}
where $m\in\mathbb{Z},$ and where $a_{j}$ are smooth functions of
$\mathbf{w},$ that is, $a_{j}\in\mathcal{D}\left(  \mathbb{S}^{n-1}\right)  .$
We call $\mathcal{D}_{\ast,\mathbf{a}}\left(  \mathbb{R}^{n}\right)  $ the
space of test functions on $\mathbb{R}^{n}$ with a thick point located at
$\mathbf{x}=\mathbf{a}.$ It is sometimes convenient to take $\mathbf{a}%
=\mathbf{0};$ we denote $\mathcal{D}_{\ast,\mathbf{0}}\left(  \mathbb{R}%
^{n}\right)  $ by $\mathcal{D}_{\ast}\left(  \mathbb{R}^{n}\right)  .$
\end{definition}

Observe that if $\phi$ is a standard test function, namely, smooth in all
$\mathbb{R}^{n}$ and with compact support, then it has a Taylor expansion,
which may be divergent, but gives a strong asymptotic expansion, $\phi\left(
\mathbf{a}+r\mathbf{w}\right)  \sim a_{0}+\sum_{j=1}^{\infty}a_{j}\left(  \mathbf{w}\right)
r^{j},$ where $a_{0}$ is just the real number $\phi\left(  \mathbf{a}\right)
.$ Hence $\mathcal{D}\left(  \mathbb{R}^{n}\right)  \subset\mathcal{D}%
_{\ast,\mathbf{a}}\left(  \mathbb{R}^{n}\right)  ;$ we denote by
\begin{equation}
i:\mathcal{D}\left(  \mathbb{R}^{n}\right)  \hookrightarrow\mathcal{D}%
_{\ast,\mathbf{a}}\left(  \mathbb{R}^{n}\right)  , \label{inclusionR^n}%
\end{equation}
the inclusion map. In fact, with the topology constructed by the following
definitions, $\mathcal{D}\left(  \mathbb{R}^{n}\right)  $ is actually a
\textsl{closed} subspace of $\mathcal{D}_{\ast,\mathbf{a}}\left(
\mathbb{R}^{n}\right)  .$

The following auxiliary spaces are needed for defining the topoloty of
$\mathcal{D}_{\ast,\mathbf{a}}\left(  \mathbb{R}^{n}\right)  $

\begin{definition}
Let $m$ be a fixed integer. The subspace $\mathcal{D}_{\ast,\mathbf{a}%
}^{\left[  m\right]  }\left(  \mathbb{R}^{n}\right)  $ consists of those test
functions $\phi$ whose expansion $\left(  \ref{test function higher dim}%
\right)  $ begins at $m.$ For a fixed compact $K$ whose interior contains
$\mathbf{a},$ $\mathcal{D}_{\ast,\mathbf{a}}^{\left[  m;K\right]  }\left(
\mathbb{R}^{n}\right)  $ is the subspace formed by those test functions of
$\mathcal{D}_{\ast,\mathbf{a}}^{\left[  m\right]  }\left(  \mathbb{R}%
^{n}\right)  $ that vanish in $\mathbb{R}^{n}\backslash K.$
\end{definition}

\begin{definition}
Let $m$ be a fixed integer and $K$ a compact subset of $\mathbb{R}^{n}$ whose
interior contains $\mathbf{a}.$ The topology of $\mathcal{D}_{\ast,\mathbf{a}%
}^{\left[  m;K\right]  }\left(  \mathbb{R}^{n}\right)  $ is given by the
seminorms $\left\{  \left\Vert {}\right\Vert _{q,s}\right\}  _{q>m,s\geq0}$
defined as
\[
\left\Vert \phi\right\Vert _{q,s}=\underset{\mathbf{x}+\mathbf{a}\in K}{\sup
}\underset{\left\vert \mathbf{p}\right\vert \leq s}{\sup}r^{-q}\left\vert \left(
\partial/\partial\mathbf{x}\right)  ^{\mathbf{p}}\phi\left(  \mathbf{a}%
+\mathbf{x}\right)  -\sum_{j=m-\left\vert \mathbf{p}\right\vert }%
^{q-1}a_{j.\mathbf{p}}\left(  \mathbf{w}\right)  r^{j}\right\vert ,
\]
where $\mathbf{x}=r\mathbf{w,p\in}\mathbb{N}^{n},$ and $\left(  \partial
/\partial\mathbf{x}\right)  ^{\mathbf{p}}\phi\left(  \mathbf{a}+\mathbf{x}%
\right)  \sim\sum_{j=m-\left\vert \mathbf{p}\right\vert }^{\infty
}a_{j.\mathbf{p}}\left(  \mathbf{w}\right)  r^{j}.$ The topology of
$\mathcal{D}_{\ast,\mathbf{a}}^{\left[  m\right]  }\left(  \mathbb{R}%
^{n}\right)  $ is the inductive limit topology of the $\mathcal{D}%
_{\ast,\mathbf{a}}^{\left[  m;K\right]  }\left(  \mathbb{R}^{n}\right)  $ as
$K\nearrow\infty.$ The topology of $\mathcal{D}_{\ast,\mathbf{a}}\left(
\mathbb{R}^{n}\right)  $ is the inductive limit topology of the $\mathcal{D}%
_{\ast,\mathbf{a}}^{\left[  m\right]  }\left(  \mathbb{R}^{n}\right)  $ as
$m\searrow-\infty.$
\end{definition}

Next we quote the definition of the "thick distribution space":

\begin{definition}
The space of distributions on $\mathbb{R}^{n}$ with a thick point at
$\mathbf{x}=\mathbf{a}$ is the dual space of $\mathcal{D}_{\ast,a}\left(
\mathbb{R}^{n}\right)  .$ We denoted it by $\mathcal{D}_{\ast,a}^{\prime}\left(
\mathbb{R}^{n}\right)  ,$ or just as $\mathcal{D}_{\ast}^{\prime}\left(
\mathbb{R}^{n}\right)  $ when $\mathbf{a}=\mathbf{0.}$ We call the elements of
$\mathcal{D}_{\ast,a}^{\prime}\left(  \mathbb{R}^{n}\right)  $ ``thick distributions".
\end{definition}

Let $\pi:\mathcal{D}_{\ast,a}^{\prime}\left(  \mathbb{R}^{n}\right)
\rightarrow\mathcal{D}^{\prime}\left(  \mathbb{R}^{n}\right)  $ denote the
projection operator, dual to the inclusion map $i:\mathcal{D}\left(
\mathbb{R}^{n}\right)  \hookrightarrow\mathcal{D}_{\ast,\mathbf{a}}\left(
\mathbb{R}^{n}\right)  .$ Since $\mathcal{D}\left(  \mathbb{R}^{n}\right)  $
is closed in $\mathcal{D}_{\ast,\mathbf{a}}\left(  \mathbb{R}^{n}\right)  ,$
the Hahn-Banach theorem immediately yields that:

\begin{theorem}
Let $f$ be any distribution in $\mathcal{D}^{\prime}\left(  \mathbb{R}%
^{n}\right)  ,$ then there exists thick distributions $g\in\mathcal{D}%
_{\ast,a}^{\prime}\left(  \mathbb{R}^{n}\right)  $ such that $\pi\left(
g\right)  =f.$
\end{theorem}

\begin{example}
(Thick delta functions of degree q) Let $g\left(  \mathbf{w}\right)  $ is a
distribution in $\mathbb{S}^{n-1}.$ The thick delta function of degree q,
denoted as $g\delta_{\ast}^{[q]},$ acts on a thick test function $\phi\left(
\mathbf{x}\right)  $ as%
\[
\left\langle g\delta_{\ast}^{[q]},\phi\right\rangle =\frac{1}{C_{n-1}%
}\left\langle g\left(  \mathbf{w}\right)  ,a_{q}\left(  \mathbf{w}\right)
\right\rangle ,
\]
where $\phi\left(  r\mathbf{w}\right)  \sim\sum_{j=m}^{\infty}a_{j}\left(
\mathbf{w}\right)  r^{j}$ as $\mathbf{x}\rightarrow0,$ and $C_{n-1}%
=\int_{\mathbb{S}^{n-1}}d\sigma\left(  \mathbf{w}\right)  $ is the surface
area of the n-dimensional unit sphere.

We consider the special case when $q=0,g\left(  \mathbf{w}\right)  \equiv1,$
we denote it $\delta_{\ast}^{[0]}:=\delta_{\ast}.$ One could easily check that
$\pi\left(  \delta_{\ast}\right)  =\delta,$ is the famous Dirac delta function.
\end{example}

All operations on thick distributions are defined via duality, especially the
multiplication and the derivative:

\begin{definition}
If $\psi\phi\in\mathcal{D}_{\ast}\left(  \mathbb{R}^{n}\right)  $ for any
$\phi\in\mathcal{D}_{\ast}\left(  \mathbb{R}^{n}\right)  ,$ then for an
$f\in\mathcal{D}_{\ast}^{\prime}\left(  \mathbb{R}^{n}\right)  ,$ $\psi f$ is
well-defined as $\left\langle \psi f,\phi\right\rangle =\left\langle
f,\psi\phi\right\rangle .$
\end{definition}

\begin{definition}
Let $f\in\mathcal{D}_{\ast}^{\prime}\left(  \mathbb{R}^{n}\right)  $ is a
thick distribution, then the thick distributional derivative of $f$ is defined
as
\[
\left\langle \frac{\partial^{\ast}f}{\partial x_{j}},\phi\right\rangle
=-\left\langle f,\frac{\partial\phi}{\partial x_{j}}\right\rangle ,\text{
\ \ \ \ \ }\phi\in\mathcal{D}_{\ast}\left(  \mathbb{R}^{n}\right)  .
\]

\end{definition}

One could easily see that $\pi\left(  \frac{\partial^{\ast}f}{\partial x_{j}%
}\right)  =\frac{\overline{\partial}f}{\partial x_{j}},$ the distributional
derivative of $f.$

\section{Reconstruction of the one-dimensional case\label{sec3}}

\subsection{Asymptotic expansions with respect to r\label{subsec3.1}}

Now let us view the line segment $[-1,1]$ as 1-dimensional unit ball, and the
boundary, i.e., the two points at $-1$ and $1$ as the \textquotedblleft%
$0\text{{}}\text{ dimensional unit sphere}$\textquotedblright. We denote the
two boundary points $\mathbf{-1}\text{ and}$ $\mathbf{1},\text{respectively.
}$We denote the set $\left\{  \mathbf{-1,1}\right\}  \text{ a}\text{s}%
\text{{}}$ $\mathbb{S}^{0}$, in accordance to the name \textquotedblleft0
dimensional unit sphere\textquotedblright.

Thus we can generalize the concept of \textquotedblleft functions on the unit
sphere\textquotedblright\ to \textquotedblleft functions on the 0 dimensional
unit sphere\textquotedblright: that is, a function from two points to
$\mathbb{R}$. Notice that the two points are disconnected, and any functions
from one point to $\mathbb{R}\text{{}}\text{ }$ is just a constant.

Now we can express the $\mathbb{R}\backslash\{0\}$ as $\mathbb{R}%
\backslash\{0\}\subset\mathbb{S}^{0}\mathbb{\times R}_{\geqslant0}:x=\left(
\mathbf{w},r\right)  ,$ where $r=\left\vert x\right\vert ;$ $\mathbf{w=1}$
when $x>0$ and $\mathbf{w=-1}$ when $x<0.$ That is, if $r>0,\left(
\mathbf{1},r\right)  $ denotes all positive numbers while $\left(
\mathbf{-1},r\right)  $ denotes all negative numbers. Notice there are two
points in $\mathbb{S}^{0}\mathbb{\times R}_{\geqslant0}$ with $r=0,$ that is,
$\left(  \mathbf{1},0\right)  $ and $\left(  \mathbf{-1},0\right)  .$

Notice that $\mathbb{S}^{0}$ has the natural discrete topology. We endow the
space $\mathbb{S}^{0}\times\mathbb{R}_{\geq0}$ with the product topology. We
endow $\mathbb{S}^{0}\times\mathbb{R}_{+}$ with the product topology, and it is
not hard to see that $\mathbb{R}\backslash\left\{  0\right\}  $ is
homeomorphic to $\mathbb{S}^{0}\times\mathbb{R}_{+}.$

\begin{definition}
\label{def1}Let $r=\left\vert x\right\vert ,\text{{}we say }$a function
$f\left(  x\right)  =f\left(  \mathbf{w},r\right)  $ defined on $\mathbb{R}%
\backslash\{0\}$ has an asymptotic expansion $\overset{\infty
}{\underset{j=m}{\sum}}a_{j}\left(  \mathbf{w}\right)  r^{j}$ , as
$x\rightarrow0\text{{}}$ , where $\mathbf{w\in}$ $\mathbb{S}^{0},$
$a_{i}\left(  \mathbf{w}\right)  $ is a function on $\mathbb{S}^{0},$ if
\begin{equation}
\underset{r\rightarrow0^{+}\text{{}}}{\lim}\left\vert f\,\left(  x\right)
-\overset{M}{\underset{j=m}{\sum}}a_{j}\left(  \mathbf{w}\right)
r^{j}\right\vert r^{-M}=0,\text{ \ uniformly on }\mathbf{w\in}\mathbb{S}^{0}.
\label{asymp1}
\end{equation}
In this case we write $f\left(  x\right)  \sim\overset{\infty
}{\underset{j=m}{\sum}}a_{j}\left(  \mathbf{w}\right)  r^{j}$ as
$x\rightarrow0.$ In fact, we can interchange $x\rightarrow0$ with
$r\rightarrow0^{+}$ here.
\end{definition}

Let us present a few examples.

\begin{example}
The first example would be polynomials. Without loss of generality, suppose
$f\left(  x\right)  =a_{2m}x^{2m}+a_{2m-1}x^{2m-1}+...+a_{0}$ is a polynomial
of even degree. One could see that $f\left(  x\right)  =f\left(
\mathbf{w},r\right)  =a_{2m}r^{2m}+a_{2m-1}\left(  \mathbf{w}\right)
r^{2m-1}+...+a_{0},$ where the coefficients of even power terms $a_{2k}$ are constants and the coefficients $a_{2k-1}\left(  \mathbf{w}%
\right)  $ of all odd power terms are functions on $\mathbb{S}^{0}:$%
\[
a_{2k-1}\left(  \mathbf{w}\right)  =\left\{
\begin{array}
[c]{c}%
a_{2k-1}\text{ \ when }\mathbf{w=1}\\
-a_{2k-1\text{ \ \ }}\text{when }\mathbf{w=-1}%
\end{array}
\right.  ,1\leq k\leq m.
\]

Clearly, we have the asymptotic expansion $f\left(  \mathbf{w},r\right)  \sim
a_{0}+...+a_{2m-1}\left(  \mathbf{w}\right)  r^{2m-1}+a_{2m}r^{2m}$ as
$x\rightarrow0.$
\end{example}

The next example shows the elegance of such description when we try to
describe functions with a jump singularity at the origin. Let us present the
Heaviside function.

\begin{example}
The Heaviside function is given as
\[
H\left(  x\right)  =\left\{
\begin{array}
[c]{c}%
1\text{ \ \ when }x>0\\
0\text{ \ \ when }x<0
\end{array}
\right. .
\]
Let us write it in the above notation, then
\begin{equation}
H\left(  x\right)  =a_{0}\left(  \mathbf{w}\right)  =\left\{
\begin{array}
[c]{c}%
1\text{ \ \ when }\mathbf{w=1}\\
0\text{ \ \ when }\mathbf{w=-1}%
\end{array}
\right.  ,\label{Heaviside1}
\end{equation}
and it admits an asymptotic expansion $H\left(  x\right)  \sim a_{0}\left(
\mathbf{w}\right)  $ as $x\rightarrow0.$
\end{example}

Let us present a lemma now concerning the relationship between the asymptotic
expansion with respect to $x$ and the asymptotic expansion with respect to
$r:$

\begin{lemma}
\label{Lemma 4.1}Suppose $f\left(  x\right)  :\mathbb{R}\backslash
\{0\}\rightarrow\mathbb{R}$ admits an asymptotic expansion $f\sim
\overset{\infty}{\underset{j=m}{\sum}}a_{2j}x^{2j}+a_{2j+1}x^{2j+1}$ as
$x\rightarrow0,$ then it admits an asymptotic expansion $f\sim\overset{\infty
}{\underset{j=m}{\sum}}a_{2j}r^{2j}+a_{2j+1}\left(  \mathbf{w}\right)
r^{2j+1}$ as $r\rightarrow0^{+}.$ Where $a_{2j}$ are constants, and 
\[a_{2j+1}\left(  \mathbf{w}\right)
=\left\{
\begin{array}
[c]{c}%
a_{2j+1}\text{ \ \ when }\mathbf{w=1}\\
-a_{2j+1}\text{ \ \ when }\mathbf{w=-1}%
\end{array}
\right.  .\]

Here we can start the expansion with an even term because we can
allow $a_{2m}=0.$

\begin{proof}
Since $f\left(  x\right)  $ admits an asymptotic expansion $f\sim
\overset{\infty}{\underset{j=m}{\sum}}a_{2j}x^{2j}+a_{2j+1}x^{2j+1},$ then
\[
\underset{x\rightarrow0}{\lim}\left\vert f\left(  x\right)
-\overset{M}{\underset{j=2m}{\sum}}a_{j}x^{j}\right\vert \left\vert
x\right\vert ^{-M}=0
\]
On the other hand,
\begin{equation}
0=\underset{x\rightarrow0^{+}}{\lim}\left\vert f\left(  x\right)
-\overset{M}{\underset{j=2m}{\sum}}a_{j}x^{j}\right\vert \left\vert
x\right\vert ^{-M}=\underset{r\rightarrow0^{+}}{\lim}\left\vert f\left(
r\right)  -\overset{M}{\underset{j=2m}{\sum}}a_{j}r^{j}\right\vert r^{-M}
\label{Lemmax}%
\end{equation}
And%
\begin{equation}
0=\underset{x\rightarrow0^{-}}{\lim}\left\vert f\left(  x\right)
-\overset{M}{\underset{j=2m}{\sum}}a_{j}x^{j}\right\vert \left\vert
x\right\vert ^{-M}=\underset{r\rightarrow0^{+}}{\lim}\left\vert f\left(
-r\right)  -\overset{M}{\underset{j=2m}{\sum}}a_{j}\left(  -r\right)
^{j}\right\vert r^{-M} \label{Lemmaxr}
\end{equation}

Combine(\ref{Lemmax}) and (\ref{Lemmaxr}), we have
\[
\underset{r\rightarrow0^{+}}{\lim}\left\vert f\left(  x\right)
-\overset{M}{\underset{j=2m}{\sum}}a_{j}\left(  \mathbf{w}\right)
r^{j}\right\vert r^{-M}=0,
\]
where $a_{j}\left(  \mathbf{w}\right)  =a_{j}$ when $j$ is even, and
\[
a_{j}\left(  \mathbf{w}\right)  =\left\{
\begin{array}
[c]{c}%
a_{j}\text{ \ \ when }\mathbf{w=1}\\
-a_{j}\text{ \ \ when }\mathbf{w=-1}%
\end{array}
\right.
\]
when $j$ is odd. Thus by definition \ref{def1}, we obtain the statement of the lemma.
\end{proof}
\end{lemma}

\begin{example}
\label{example3} Let $f:\mathbb{R}
\rightarrow\mathbb{R}$ be a smooth function. It has a Taylor expansion at the
origin:\ $f\left(  x\right)  \sim\overset{\infty}{\underset{j=0}{\sum}}%
\frac{f^{\left(  j\right)  }\left(  0\right)  }{j!}x^{j}.$ It may not
converge. But it is an asymptotic expansion as $x\rightarrow0.$ By the Lemma
\ref{Lemma 4.1}, we have an asymptotic expansion
\[
f\sim\overset{\infty}{\underset{j=0}{\sum}}a_{2j}r^{2j}+a_{2j+1}\left(
\mathbf{w}\right)  r^{2j+1}\text{ \ \ as }r\rightarrow0^{+},
\]
where $a_{2j}=\frac{f^{\left(  2j\right)  }\left(  0\right)  }{\left(
2j\right)  !}$ and
\[
a_{2j+1}\left(  \mathbf{w}\right)  =\left\{
\begin{array}
[c]{c}%
\frac{f^{\left(  2j+1\right)  }\left(  0\right)  }{\left(  2j+1\right)
!}\text{ \ \ when }\mathbf{w=1}\\
-\frac{f^{\left(  2j+1\right)  }\left(  0\right)  }{\left(  2j+1\right)
!}\text{ \ \ when }\mathbf{w=-1}%
\end{array}
\right. .
\]

\end{example}

\bigskip

Next we want to discuss the derivative of $\overset{\infty
}{\underset{j=m}{\sum}}a_{j}\left(  \mathbf{w}\right)  r^{j}$ with respect to
$x.$ We first discuss $\frac{d\left(  a\left(  \mathbf{w}\right)
r^{j}\right)  }{dx},$ where $\mathbf{w\in}$ $\mathbb{S}^{0},$ $a\left(
\mathbf{w}\right)  $ is a function on $\mathbb{S}^{0}.$ Because there is a
natural inclusion $\mathbb{R}\backslash\{0\}\subset\mathbb{S}^{0}%
\mathbb{\times R}_{\geqslant0},$ clearly, $a\left(  \mathbf{w}\right)  r^{j}$
can be viewed as a function on $x$ when $r=\left\vert x\right\vert \neq0.$
Thus, it is legal to talk about "derivative with respect to x" in the usual sense.

We first discuss the derivative $\frac{d\left(  a\left(  \mathbf{w}\right)
r^{j}\right)  }{dx},j\neq0,$ at $x_{0}>0.$ Denote the coordinate of $x_{0}$ at
$\mathbb{S}^{0}\times\mathbb{R}_{\geq0}$ as $\left(  \mathbf{1},x_{0}\right)
.$ Since $\mathbb{S}^{0}$ is endowed with the discrete topology, there is a
small neighborhood of $\left(  \mathbf{1},x_{0}\right) ,$ denoted as $\left\{
\mathbf{1}\right\}  \times\lbrack x_{0}-\delta,x_{0}+\delta],$ on which the
function $a\left(  \mathbf{w}\right)  r^{j}$ equals $a\left(  \mathbf{1}%
\right)  x^{j},$ where $a\left(  \mathbf{1}\right)  $ is the value of the
function $a\left(  \mathbf{w}\right)  $ at $\mathbf{w}=\mathbf{1}$, namely, a
constant. Thus we have
\[
\left.  \frac{d\left(  a\left(  \mathbf{w}\right)  r^{j}\right)  }%
{dx}\right\vert _{x=x_{0}}=\left.  \frac{d\left(  a\left(  \mathbf{1}\right)
x^{j}\right)  }{dx}\right\vert _{x=x_{0}}=\left.  a\left(  \mathbf{1}\right)
jx^{j-1}\right\vert _{x=x_{0}}.
\]
Similarly, when $x_{0}<0,$ there is a small neighborhood of $\left(
-\mathbf{1},-x_{0}\right) ,$ denoted as $\left\{  -\mathbf{1}\right\}
\times\lbrack-x_{0}-\delta,-x_{0}+\delta],$ on which the function $a\left(
\mathbf{w}\right)  r^{j}$ equals $a\left(  -\mathbf{1}\right)  \left(
-x\right)  ^{j},$ hence when $x_{0}<0,$ $j\neq0,$ we have
\[
\left.  \frac{d\left(  a\left(  \mathbf{w}\right)  r^{j}\right)  }%
{dx}\right\vert _{x=x_{0}}=\left.  -a\left(  -\mathbf{1}\right)  j\left(
-x\right)  ^{j-1}\right\vert _{x=x_{0}}.
\]
When $j=0,$ a similar analysis shows that $\left.  \frac{d\left(  a\left(
\mathbf{w}\right)  \right)  }{dx}\right\vert _{x=x_{0}}=0$ when $x_{0}\neq0.$

When $x=0,$ $a\left(  \mathbf{w}\right)  $ can be viewed as a multi-valued
function on $x.$ On the other hand, if $j>0,$ then $a\left(  \mathbf{w}%
\right)  r^{j}=0$ at $x=0.$ We can talk about the so called "left-derivative"
and "right-derivative" of $a\left(  \mathbf{w}\right)  r^{j}$, $j>0$ at $x=0$.
It is not hard to see, the "left-derivative" and "right-derivative" are both
$0$ when $j>1;$ and the "right-derivative" equals $a\left(  \mathbf{1}\right)
$ while the "left-derivative" equals $a\left(  -\mathbf{1}\right)  $ when
$j=1,$ at $x=0.$

Based on the above discussion, we have the following theorem regarding the
derivative $\frac{d\left(  a\left(  \mathbf{w}\right)  r^{j}\right)  }{dx}:$

\begin{theorem}
Let $a\left(  \mathbf{w}\right)  r^{j}$ be a function on $\mathbb{S}^{0}%
\times\mathbb{R}_{+},$ then
\[
\frac{d\left(  a\left(  \mathbf{w}\right)  r^{j}\right)  }{dx}=\left\{
\begin{array}
[c]{c}
a\left(  \mathbf{1}\right)  jx^{j-1}\text{ \ \ \ \ \ \ when }x>0\\
-a\left(  -\mathbf{1}\right)  j\left(  -x\right)  ^{j-1}\text{ \ when }x<0
\end{array}
\right. .
\]
In summary, it amounts to saying $\frac{d\left(  a\left(  \mathbf{w}\right)
r^{j}\right)  }{dx}=b\left(  \mathbf{w}\right)  r^{j-1},$ where $b\left(
\mathbf{w}\right)  $ is a function on $\mathbb{S}^{0}:$
\[
b\left(  \mathbf{w}\right)  =\left\{
\begin{array}
[c]{c}
a\left(  \mathbf{1}\right)  j\text{ \ \ \ \ \ \ when }\mathbf{w}=\mathbf{1}\\
-a\left(  -\mathbf{1}\right)  j\text{ \ \ \ \ \ when }\mathbf{w}=-\mathbf{1}%
\end{array}
\right. .
\]

\end{theorem}

We immediately have the following theorem:

\begin{theorem}
If $f\left(  x\right)  \sim\overset{\infty}{\underset{j=m}{\sum}}a_{j}\left(
\mathbf{w}\right)  r^{j}$ as $x\rightarrow0,$ then the term-by-term derivative
with respect to $x$ of the expasion takes the following form in $\mathbb{S}%
^{0}\times\mathbb{R}_{+}:$
\begin{align}
\overset{\infty}{\underset{j=m}{\sum}}\frac{d\left(  a_{j}\left(
\mathbf{w}\right)  r^{j}\right)  }{dx}  &  =\overset{\infty
}{\underset{j=m-1}{\sum}}a_{j,1}\left(  \mathbf{w}\right)  r^{j},\text{
}\label{derivative1}\\
\text{where }a_{j,1}\left(  \mathbf{w}\right)   &  =\left\{
\begin{array}
[c]{c}%
a_{j+1}\left(  \mathbf{w}\right)  \left(  j+1\right)  \text{ \ \ \ \ \ \ when
}\mathbf{w}=\mathbf{1}\\
-a_{j+1}\left(  \mathbf{w}\right)  \left(  j+1\right)  \text{ \ \ \ \ \ when
}\mathbf{w}=-\mathbf{1}%
\end{array}
\right.  .
\end{align}

\end{theorem}

\begin{example}
The above theorems show that the usual derivative (NOT\ the distributional
derivative) of the Heaviside function is $0$ on $\mathbb{S}^{0}\times
\mathbb{R}_{+}=\mathbb{R}\backslash\left\{  0\right\}  .$ Here we can clearly
distinguish between the "usual derivative" and the "distributional
derivative". Since we know the famous fact that the "distributinal derivative"
of the Heaviside function is the Dirac delta function.
\end{example}

Now we are ready to construct the space of thick test functions.

\subsection{Reconstruction of the space of test functions on $\mathbb{R}$ with
a thick point\label{subsec3.2}}

\begin{definition}
\label{strong}Let $\phi\in C^{\infty}\left(  \mathbb{R}\backslash\left\{
0\right\}  \right)  .$ We say that $\phi\left(  x\right)$ has a ``strong" expansion $
\sim\overset{\infty}{\underset{j=m}{\sum}}a_{j}\left(  \mathbf{w}\right)
r^{j}$ as $x\rightarrow0$ if the asymptotic development of
$\frac{d^{p}\phi}{dx^{p}}$ exists and equals the term-by-term differentiation
$\overset{\infty}{\underset{j=m}{\sum}}\frac{d^p(a_{j}\left(  \mathbf{w}\right)
r^{j})}{dx^p},$ as the derivative was discussed in the previous section.
\end{definition}

We now define the space of test functions with a thick point at $x=a.$

\begin{definition}
\label{test functions}Let $\mathcal{D}_{\ast,a}\left(  \mathbb{R}\right)  $
denote the vector space of all compactly supported smooth functions $\phi$ defined in $\mathbb{R}\backslash\left\{  a\right\}  ,$ that admit a strong asymptotic expansion of the form
\begin{equation}
\phi\left(  a+x\right)  \sim\underset{j=m}{\overset{\infty}{\sum}}a_{j}\left(
\mathbf{w}\right)  r^{j},\text{ \ \ \ \ \ as }x\rightarrow0.\label{asymp2}%
\end{equation}
\bigskip where $a_{j}\left(  \mathbf{w}\right)  $ is a function on
$\mathbb{S}^{0}$ as defined above in definition \ref{def1}. We call
$\mathcal{D}_{\ast,a}\left(  \mathbb{R}\right)  $ ``the space of test functions
on $\mathbb{R}$ with a thick point located at $x=a".$ We denote $\mathcal{D}_{\ast,0}\left(  \mathbb{R}\right)  $ as $\mathcal{D}_{\ast}\left(
\mathbb{R}\right)  .$
\end{definition}

The definition of $\mathcal{D}_{\ast,a}\left(  \mathbb{R}\right)  $ is very
much analogous to the defintion of $\mathcal{D}_{\ast,a}\left(  \mathbb{R}%
^{n}\right)  $ when $n\geq2.$ \cite{Thick}

On the other hand, the thick test functions defined in \cite{Thick1} is
included in the above definition. Recall that a thick test function defined in
\cite{Thick1} is a compactly supported function $\phi$ with domain
$\mathbb{R},$ smooth in $\mathbb{R}\backslash\left\{  a\right\}  ,$ and at
$x=a$ all its one-sided derivatives,
\[
\phi^{\left(  n\right)  }\left(  a\pm0\right)  =\underset{x\rightarrow a^{\pm
}}{\lim}\phi^{\left(  n\right)  }\left(  x\right)  ,\text{ \ \ }\forall
n\in\mathbb{N},
\]
exist. Here let us introduce a different notation $\mathcal{D}_{\ast,a}%
^{old}\left(  \mathbb{R}\right)  $ to denote the space of such functions. One
can see that any function in $\mathcal{D}_{\ast,a}^{old}\left(  \mathbb{R}%
\right)  $ admits a strong asymptotic expansion
\begin{align}
\phi\left(  a+x\right)   &  \sim\underset{j=0}{\overset{\infty}{\sum}}%
a_{j}\left(  \mathbf{w}\right)  r^{j},\text{ \ \ \ as }x\rightarrow
0,\label{thick1def}\\
where\text{ }a_{j}\left(  \mathbf{w}\right)   &  =\left\{
\begin{array}
[c]{c}%
\frac{\phi^{\left(  j\right)  }\left(  a+0\right)  }{j!}\text{ \ \ \ \ when
}\mathbf{w=1}\\
\left(  -1\right)  ^{j}\frac{\phi^{\left(  j\right)  }\left(  a-0\right)
}{j!}\text{ \ \ \ \ when }\mathbf{w=-1}%
\end{array}
\right.  .\nonumber
\end{align}
In particular, if $\phi\left(  x\right)  $ has a jump discontinuity at $x=a,$
then in the expansion (\ref{thick1def}), $a_{0}\left(  \mathbf{1}\right)  \neq
a_{0}\left(  \mathbf{-1}\right)  .$

By example \ref{example3}, all compactly supported smooth functions form a
closed subspace of $\mathcal{D}_{\ast,a}^{old}\left(  \mathbb{R}\right)  ,$
thus there are natural inclusion maps $:$%
\begin{equation}
\mathcal{D}\left(  \mathbb{R}\right)  \hookrightarrow\mathcal{D}_{\ast
,a}^{old}\left(  \mathbb{R}\right)  \hookrightarrow\mathcal{D}_{\ast,a}\left(
\mathbb{R}\right)  . \label{inclusion}%
\end{equation}

\bigskip

Now let us define a topology on $\mathcal{D}_{\ast,a}\left(  \mathbb{R}%
\right)  $ and make it a topological vector space.

Similar to \cite{Thick}, we denote the subspace $\mathcal{D}_{\ast,a}^{\left[
m\right]  }\left(  \mathbb{R}\right)  $ as the test functions $\phi$ whose
expansion (\ref{asymp2}) begins at $m.$ For a fixed compact $K$ whose interior
contains $a,\mathcal{D}_{\ast,a}^{\left[  m,K\right]  }\left(  \mathbb{R}%
\right)  $ is the subspace formed by those test functions of $\mathcal{D}%
_{\ast,a}^{\left[  m\right]  }\left(  \mathbb{R}\right)  $ that vanish in
$\mathbb{R}\backslash K.$

Now let us define the topology of the space of thick test functions in
$\mathbb{R}.$

\begin{definition}
\label{Deftopo}Let $m$ be a fixed integer and $K$ a compact subset of
$\mathbb{R}$ whose interior contains $a.$ The topology of $\mathcal{D}%
_{\ast,a}^{\left[  m,K\right]  }\left(  \mathbb{R}\right)  $ is given by the
seminorms $\left\{  \left\Vert {}\right\Vert _{q,s}\right\}  _{s\geq0}$
defined as
\[
\left\Vert \phi\right\Vert _{q,s}=\underset{x+a\in K}{\sup}\underset{0\leq
p\leq s}{\sup}r^{-q}\left\vert \left(  d/dx\right)  ^{p}\phi\left(
a+x\right)  -\overset{q-1}{\underset{j=m-p}{\sum}}a_{j,p}\left(
\mathbf{w}\right)  r^{j}\right\vert ,
\]
where $\left(  d/dx\right)  ^{p}\phi\left(  a+x\right)  \sim
\overset{q-1}{\underset{j=m-p}{\sum}}a_{j,p}\left(  \mathbf{w}\right)  r^{j}.$
The topology of $\mathcal{D}_{\ast,a}^{\left[  m\right]  }\left(
\mathbb{R}\right)  $ is the inductive limit topology of the $\mathcal{D}%
_{\ast,a}^{\left[  m,K\right]  }\left(  \mathbb{R}\right)  $ as $K\nearrow
\infty.$ The topology of $\mathcal{D}_{\ast,a}\left(  \mathbb{R}\right)  $ is
the inductive limit topology of $\mathcal{D}_{\ast,a}^{\left[  m\right]
}\left(  \mathbb{R}\right)  $ as $m\searrow-\infty.$
\end{definition}


One could see that the ``space of test functions in $\mathbb{R}$" that is
introduced in the paper \cite{Thick1} is a closed subspace of
$\mathcal{D}_{\ast,a}\left(  \mathbb{R}\right)  $ introduced above:
$\mathcal{D}_{\ast,a}^{old}\left(  \mathbb{R}\right)  \subseteq\mathcal{D}%
_{\ast,a}\left(  \mathbb{R}\right)  .$ Moreover, $\mathcal{D}_{\ast,a}%
^{old}\left(  \mathbb{R}\right)  $ is closed in $\mathcal{D}_{\ast,a}\left(
\mathbb{R}\right)  $ with respect to derivatives.

Furthermore, there is a natural inclusion map $i:\mathcal{D}\left(
\mathbb{R}\right)  \hookrightarrow\mathcal{D}_{\ast,a}\left(  \mathbb{R}%
\right)  $. With the topology defined in definition \ref{Deftopo}, the space
of usual test functions $\mathcal{D}\left(  \mathbb{R}\right)  $ is a closed
subspace of $\mathcal{D}_{\ast,a}\left(  \mathbb{R}\right)  ,$ similar to the
case of higher dimensions \cite{Thick}.

\bigskip

\subsection{Space of distributions on $\mathbb{R}$ with a thick
point\label{subsec3.3}}

\bigskip

Having all the preparation above, we can define distributions on $\mathbb{R}$
with a thick point.

\begin{definition}
The space of distributions on $\mathbb{R}$ with a thick point at $x=a$ is the
dual space of $\mathcal{D}_{\ast,a}\left(  \mathbb{R}\right)  .$ We denote it
$\mathcal{D}_{\ast,a}^{\prime}\left(  \mathbb{R}\right)  ,$ or just as
$\mathcal{D}_{\ast}^{\prime}\left(  \mathbb{R}\right)  $ when $a=0.$ We call
the elements of $\mathcal{D}_{\ast,a}^{\prime}\left(  \mathbb{R}\right)  $
``thick distributions".
\end{definition}

Let $\pi:\mathcal{D}_{\ast,a}^{\prime}\left(  \mathbb{R}\right)
\rightarrow\mathcal{D}^{\prime}\left(  \mathbb{R}\right)  ,$ be the projection
operator, dual to the inclusion $i:\mathcal{D}\left(  \mathbb{R}\right)
\hookrightarrow\mathcal{D}_{\ast,a}\left(  \mathbb{R}\right)  .$ Since
$\mathcal{D}\left(  \mathbb{R}\right)  $ is closed in $\mathcal{D}_{\ast
,a}\left(  \mathbb{R}\right)  ,$ by the Hanh-Banach theorem we have the following result.

\begin{theorem}
Let $f$ be any distribution in $\mathcal{D}^{\prime}\left(  \mathbb{R}\right)
,$ then there exist thick distributions $g\in\mathcal{D}_{\ast,a}^{\prime
}\left(  \mathbb{R}\right)  $ such that $\pi\left(  g\right)  =f.$
\end{theorem}

Naturally, if $f\in\mathcal{D}^{\prime}\left(  \mathbb{R}\right)  $ then there
are infinitely many thick distributions \ $g$ with $\pi\left(  g\right)  =f.$ 

Before giving any examples of thick distributions, let us recall that by the
convention of the discrete measure, the "integral" of a function on a discrete
set is just the summation of the function over these discrete points. Using
our notation in this present article, one just write, for any function $\phi$
defined on $\mathbb{S}^{0}$:%
\[
\int_{S^{0}}\phi\left(  \mathbf{w}\right)  d\sigma\left(  \mathbf{w}\right)
=\phi\left(  \mathbf{1}\right)  +\phi\left(  \mathbf{-1}\right)  .
\]
For example, for the Heaviside function as in equation (\ref{Heaviside1}):
$H\left(  x\right)  =a_{0}\left(  \mathbf{w}\right)  ,$%
\begin{equation}
\int_{S^{0}}a_{0}\left(  \mathbf{w}\right)  d\sigma\left(  \mathbf{w}\right)
=1. \label{integralHeaviside}%
\end{equation}
For a constant function $\phi\left(  \mathbf{w}\right)  \equiv1,$%
\begin{equation}
\int_{S^{0}}1d\sigma\left(  \mathbf{w}\right)  =\int_{S^{0}}d\sigma\left(
\mathbf{w}\right)  =2. \label{unitmeasure}%
\end{equation}

Thus one can discuss the "double integral"\ on $\mathbb{S}^{0}\times
\mathbb{R}_{+}$ if it exists$:$%
\begin{align}
\int_{S^{0}}\int_{0}^{+\infty}\phi\left(  \mathbf{w},r\right)  drd\sigma
\left(  \mathbf{w}\right)   &  =\int_{0}^{+\infty}\left[  \phi\left(
\mathbf{1},r\right)  +\phi\left(  \mathbf{-1},r\right)  \right]
dr\label{doubleintegral}\\
&  =\int_{0}^{+\infty}\phi\left(  x\right)  dx+\int_{0}^{-\infty}\phi\left(
x\right)  d\left(  -x\right) \nonumber\\
&  =\int_{0}^{+\infty}\phi\left(  x\right)  dx+\int_{-\infty}^{0}\phi\left(
x\right)  dx.\nonumber
\end{align}
Clearly, if $\phi\in\mathcal{D}\left(  \mathbb{R}\right)  $ is a usual test
function, then
\[
\int_{S^{0}}\int_{0}^{+\infty}\phi\left(  \mathbf{w},r\right)  drd\sigma
\left(  \mathbf{w}\right)  =\int_{-\infty}^{+\infty}\phi\left(  x\right)  dx,
\]
is just a normal integral over $\mathbb{R}$.

\bigskip

Next let us present a few examples of thick distributions in $\mathcal{D}
_{\ast,a}^{\prime}\left(  \mathbb{R}\right)  .$

It is well known that any locally integrable function $f$ defined in
$\mathbb{R}$ yields a distribution, usually denoted by the same notation $f,$
by the prescription%
\begin{equation}
\left\langle f,\phi\right\rangle =\int_{-\infty}^{+\infty}f\left(  x\right)
\phi\left(  x\right)  dx,\text{ \ \ }\phi\in\mathcal{D}\left(  R\right)  .
\label{locallyint}%
\end{equation}
Similar to the higher dimensional case, if $a\notin$ supp$f,$ that is, if
$f\left(  x\right)  =0$ for $\left\vert x-a\right\vert <\epsilon$ for some
$\epsilon>0,$ then (\ref{locallyint}) will also work in $\mathcal{D}_{\ast
,a}^{\prime}\left(  \mathbb{R}\right)  .$ However, if $a\in$ supp$f,$ then, in
general the integral $\int_{-\infty}^{+\infty}f\left(  x\right)  \phi\left(
x\right)  dx$ would be divergent and thus a thick distribution that one could
call \textquotedblleft$f$\textquotedblright\ cannot be defined in a canonical
way. Nevertheless, it is possible in many cases to define a "finite part"
distribution $Pf\left(  f\left(  x\right)  \right)  $ which is the canonical
thick distribution corresponding to $f.$ Let us recall at this point the
definition of the finite part of a limit \cite{Greenbook}.

\begin{definition}
Let $F$ be a function defined in an interval of the form $\left(  0,a\right)
$ for some $a>0.$ We say that the finite part of the limit of $F\left(
\epsilon\right)  $ as $\epsilon\rightarrow0^{+}$ exists and equals $A,$ and
denote this as $F.p.\lim_{\epsilon\rightarrow0^{+}}F\left(  \epsilon\right)
=A,$ if $F$ has the decomposition $F\left(  \epsilon\right)  =F_{fin}\left(
\epsilon\right)  +F_{\text{infin}}\left(  \epsilon\right)  ,$ where the
infinite part $F_{\text{infin}}\left(  \epsilon\right)  $ is a finite linear
combination of functions of the type $\epsilon^{-p}\ln^{q}\epsilon,$ where
$p\geq0$ and $q>0,$ and where the finite part $F_{fin}\left(  \epsilon\right)
$ is a function whose limit as $\epsilon\rightarrow0^{+}$ equals $A.$
\end{definition}

\begin{definition}
Let $f$ be a locally integrable function definied in $\mathbb{R}%
\backslash\left\{  a\right\}  .$ The thick distribution $Pf\left(  f\left(
x\right)  \right)  $ is defined as
\begin{align}
\left\langle Pf\left(  f\left(  x\right)  \right)  ,\phi\left(  x\right)
\right\rangle  &  =F.p.\int_{-\infty}^{+\infty}f\left(  x\right)  \phi\left(
x\right)  dx\label{finitepart}\\
&  =F.p.\underset{\epsilon\rightarrow0^{+}}{\lim}\int_{\left\vert
x-a\right\vert \geq\epsilon}f\left(  x\right)  \phi\left(  x\right)  dx,\text{
\ \ }\phi\in\mathcal{D}_{\ast,a}\left(  \mathbb{R}\right)  ,\nonumber
\end{align}
provided that the finite part integrals exist for all $\phi\in\mathcal{D}%
_{\ast,a}\left(  \mathbb{R}\right)  .$ Here if we set $x-a=y,\left\vert
y\right\vert =r,$ then
\begin{align}
\int_{\left\vert x-a\right\vert \geq\epsilon}f\left(  x\right)  \phi\left(
x\right)  dx  &  =\int_{\left\vert y\right\vert \geq\epsilon}f\left(
y+a\right)  \phi\left(  y+a\right)  dy=\int_{\left\vert y\right\vert
\geq\epsilon}g\left(  y\right)  \psi\left(  y\right)  dy\label{finitepart2}\\
&  =\int_{S^{0}}\int_{\epsilon}^{+\infty}g\left(  \mathbf{w},r\right)
\psi\left(  \mathbf{w},r\right)  dr\nonumber\\
&  =\int_{\epsilon}^{+\infty}g\left(  \mathbf{1},r\right)  \psi\left(
\mathbf{1},r\right)  dr+\int_{\epsilon}^{+\infty}g\left(  \mathbf{-1}%
,r\right)  \psi\left(  \mathbf{-1},r\right)  dr\nonumber\\
&  =\int_{\epsilon}^{+\infty}g\left(  y\right)  \psi\left(  y\right)
dy+\int_{-\epsilon}^{-\infty}g\left(  y\right)  \psi\left(  y\right)  d\left(
-y\right) \nonumber\\
&  =\int_{a+\epsilon}^{+\infty}f\left(  x\right)  \phi\left(  x\right)
dx+\int_{-\infty}^{a-\epsilon}f\left(  x\right)  \phi\left(  x\right)
dx.\nonumber
\end{align}

\end{definition}

\bigskip

Similar to the higher dimensional case, although the finite part limit is not
defined for all locally integrable functions $f,$ $Pf\left(  f\left(
x\right)  \right)  $ is defined in many important and interesting cases.

\begin{example}
If $\lambda\in\mathbb{C}$ then $Pf\left(  \left\vert x-a\right\vert ^{\lambda
}\right)  $ is a well-defined thick distribution of $\mathcal{D}_{\ast
,a}^{\prime}\left(  \mathbb{R}\right)  .$ Indeed, one needs to consider the
finite part of the integral $\int_{-\infty}^{+\infty}\left\vert x-a\right\vert
^{\lambda}\phi\left(  x\right)  dx$ for any $\phi\left(  x\right)
\in\mathcal{D}_{\ast,a}\left(  \mathbb{R}\right)  .$ Explicitly, since
$F.p.\int_{0}^{A}r^{\alpha}dr=A^{\alpha+1}/\left(  \alpha+1\right)
,\alpha\neq-1,F.p.\int_{0}^{A}r^{-1}dr=\log A,$ we obtain that if
$\lambda\notin\mathbb{Z}$ then%
\begin{align}
\left\langle Pf\left(  \left\vert x-a\right\vert ^{\lambda}\right)
,\phi\left(  x\right)  \right\rangle  &  =\int_{\left\vert x-a\right\vert \geq
A}\left\vert x-a\right\vert ^{\lambda}\phi\left(  x\right)  dx\label{r^lambda}%
\\
&  +\int_{\left\vert x-a\right\vert <A}\left\vert x-a\right\vert ^{\lambda
}\left(  \phi\left(  x\right)  -\underset{j\leq-\operatorname{Re}%
\lambda-1}{\sum}a_{j}\left(  w\right)  \left\vert x-a\right\vert ^{j}\right)
dx\nonumber\\
&  +\underset{j\leq-\operatorname{Re}\lambda-1}{\sum}\left(  a_{j}\left(
\mathbf{1}\right)  +a_{j}\left(  \mathbf{-1}\right)  \right)  \frac
{A^{\lambda+j+1}}{\lambda+j+1},\nonumber
\end{align}
while if $\lambda=k\in\mathbb{Z}$ then
\begin{align}
\left\langle Pf\left(  \left\vert x-a\right\vert ^{\lambda}\right)
,\phi\left(  x\right)  \right\rangle  &  =\int_{\left\vert x-a\right\vert \geq
A}\left\vert x-a\right\vert ^{k}\phi\left(  x\right)  dx\label{r^k}\\
&  +\int_{\left\vert x-a\right\vert <A}\left\vert x-a\right\vert ^{k}\left(
\phi\left(  x\right)  -\underset{j\leq-k-1}{\sum}a_{j}\left(  w\right)
\left\vert x-a\right\vert ^{j}\right)  dx\nonumber\\
&  +\underset{j<-k-1}{\sum}\left(  a_{j}\left(  \mathbf{1}\right)
+a_{j}\left(  \mathbf{-1}\right)  \right)  \frac{A^{\lambda+j+1}}{\lambda
+j+1}+\left(  a_{-k-1}\left(  \mathbf{1}\right)  +a_{-k-1}\left(
\mathbf{-1}\right)  \right)  \log A.\nonumber
\end{align}
Formulas $\left(  \ref{r^lambda}\right)  $ and (\ref{r^k}) hold for any $A>0.$
The finite part is needed for all $\lambda$ in the space of thick
distributions $\mathcal{D}_{\ast,a}^{\prime}\left(  \mathbb{R}\right)  .$
\end{example}

Using the ideas of the previous example one can show that when $\psi$ is
smooth in all of $\mathbb{R},$ i.e., $\psi\in\mathcal{E}\left(  \mathbb{R}%
\right)  ,$ then $Pf\left(  \psi\left(  x\right)  \right)  \in\mathcal{D}%
_{\ast,a}^{\prime}\left(  \mathbb{R}\right)  ;$ notice that the finite part
regularization is always needed if there is a thick point in the support of
$\psi.$

\begin{example}
\label{HeavisideFunction}This example will be the ``finite part regularization"
of the Heaviside function. The regulariztion is needed because of the
singularity of the thick test functions: let $\phi\left(  x\right)
\in\mathcal{D}_{\ast}\left(  \mathbb{R}\right)  ,$ the $Pf\left(  H\left(
x\right)  \right)  \in\mathcal{D}_{\ast}^{\prime}\left(  \mathbb{R}\right)  $
is defined as:%
\begin{align}
\left\langle Pf\left(  H\left(  x\right)  \right)  ,\phi\left(  x\right)
\right\rangle  &  =\int_{\left\vert x-a\right\vert \geq A}H\left(  x\right)
\phi\left(  x\right)  dx+\int_{\left\vert x-a\right\vert <A}H\left(  x\right)
\left(  \phi\left(  x\right)  -\underset{j\leq-1}{\sum}a_{j}\left(  w\right)
r^{j}\right)  dx\label{Heaviside}\\
&  +\underset{j<-1}{\sum}a_{j}\left(  \mathbf{1}\right)  \frac{A^{j+1}}%
{j+1}+a_{-1}\left(  \mathbf{1}\right)  \log A.\nonumber
\end{align}

\end{example}

\bigskip

Similar to the higher dimensional case, we can also define a ``thick delta
function with degree q". This definition is closely related to the extension
of Dirac delta function introduce in \cite{Thick1}. Note that here we use the
notation $C_{0}=\int_{S^{0}}1d\sigma\left(  \mathbf{w}\right)  =\int_{S^{0}%
}d\sigma\left(  \mathbf{w}\right)  =2.$

\begin{definition}
\label{thick delta function of degree q} Let $g\left(  \mathbf{w}\right)  $ be
a distribution in $\mathbb{S}^{0},$ the thick delta function of degree q,
denoted as $g\delta_{\ast}^{\left[  q\right]  },$ acts on a thick test
function $\phi\left(  x\right)  $ as
\[
\left\langle g\delta_{\ast}^{\left[  q\right]  },\phi\right\rangle
_{\mathcal{D}_{\ast}^{\prime}\left(  \mathbb{R}\right)  \times\mathcal{D}%
_{\ast}\left(  \mathbb{R}\right)  }=\frac{1}{C_{0}}\left\langle g\left(
\mathbf{w}\right)  ,a_{q}\left(  \mathbf{w}\right)  \right\rangle
_{\mathcal{D}_{\ast}^{\prime}\left(  \mathbb{S}\right)  \times\mathcal{D}%
_{\ast}\left(  \mathbb{S}\right)  },
\]
where $\phi\left(  x\right)  \sim\underset{j=m}{\overset{\infty}{\sum}}%
a_{j}\left(  \mathbf{w}\right)  r^{j},$as $x\rightarrow0,$ and $C_{0}=2$.

The thick delta function of degree 0, namely $g\delta_{\ast}^{\left[
0\right]  }$ will be denoted as $g\delta_{\ast},$ as $g\left(  \mathbf{w}%
\right)  \delta_{\ast},$ or as $g\left(  \mathbf{w}\right)  \delta_{\ast
}\left(  x\right)  .$ In particular, if $g\left(  x\right)  \equiv1,$ then we
obtain the one-dimensional ``plain thick delta function" $\delta_{\ast},$ given
as%
\[
\left\langle \delta_{\ast},\phi\right\rangle _{\mathcal{D}_{\ast}^{\prime
}\left(  \mathbb{R}\right)  \times\mathcal{D}_{\ast}\left(  \mathbb{R}\right)
}=\frac{1}{C_{0}}\int_{\mathbb{S}^{0}}a_{0}\left(  \mathbf{w}\right)
d\sigma\left(  \mathbf{w}\right)  =\frac{a_{0}\left(  \mathbf{1}\right)  }%
{2}+\frac{a_{0}\left(  \mathbf{-1}\right)  }{2}.
\]

\end{definition}

\begin{remark}
\label{remark_delta_projection} If $\phi\in\mathcal{D}\left(  \mathbb{R}%
\right)  $ is a usual test function, then by example \ref{example3} we see
that
\[
\left\langle \pi\left(  \delta_{\ast}\right)  ,\phi\right\rangle =\left\langle
\delta_{\ast},i\left(  \phi\right)  \right\rangle =\frac{\phi\left(  0\right)
}{2}+\frac{\phi\left(  0\right)  }{2}=\phi\left(  0\right)  ,
\]
hence $\pi\left(  \delta_{\ast}\right)  =\delta.$

Moreover, if $\phi\in\mathcal{D}_{\ast}^{old}\left(  \mathbb{R}\right)  ,$ let
$\pi^{\prime}:\mathcal{D}_{\ast}^{\prime}\left(  \mathbb{R}\right)
\rightarrow\mathcal{D}_{\ast}^{\prime old}\left(  \mathbb{R}\right)  $ be the
projection onto the space of thick distributions defined in \cite{Thick1}, let
$\phi_{+}\left(  0\right)  $ denote $\underset{x\rightarrow0^{+}}{\lim}%
\phi\left(  x\right)  $ and $\phi_{-}\left(  0\right)  $ denote
$\underset{x\rightarrow0^{-}}{\lim}\phi\left(  x\right)  .$ Then by equation
(\ref{thick1def}), one can see that
\begin{align*}
\left\langle \pi^{\prime}\left(  \delta_{\ast}\right)  ,\phi\right\rangle  &
=\frac{1}{2}\phi_{+}\left(  0\right)  +\frac{1}{2}\phi_{-}\left(  0\right) \\
&  =\left\langle \widetilde{\delta}\left(  x\right)  ,\phi\left(  x\right)
\right\rangle ,
\end{align*}
where $\widetilde{\delta}\left(  x\right)  $ is defined in \cite{Thick1},
equation (3.10).
\end{remark}

\begin{example}
\label{Example of thick delta}Let $g_{\lambda}\left(  \mathbf{w}\right)  $ be
a distribution in $\mathbb{S}^{0}:\left\langle g_{\lambda}\left(
\mathbf{w}\right)  ,a\left(  \mathbf{w}\right)  \right\rangle =2\lambda
a\left(  \mathbf{1}\right)  +2\left(  1-\lambda\right)  a\left(
\mathbf{-1}\right)  $, where $0\leq\lambda\leq1$ is a constant. Then
\[
\left\langle g_{\lambda}\delta_{\ast}^{\left[  q\right]  },\phi\right\rangle
_{\mathcal{D}_{\ast}^{\prime}\left(  \mathbb{R}\right)  \times\mathcal{D}%
_{\ast}\left(  \mathbb{R}\right)  }=\lambda a_{q}\left(  \mathbf{1}\right)
+\left(  1-\lambda\right)  a_{q}\left(  \mathbf{-1}\right)  .
\]
In particular, if $\lambda=1,$%
\[
\left\langle g_{1}\delta_{\ast}^{\left[  q\right]  },\phi\right\rangle
_{\mathcal{D}_{\ast}^{\prime}\left(  \mathbb{R}\right)  \times\mathcal{D}%
_{\ast}\left(  \mathbb{R}\right)  }=a_{q}\left(  \mathbf{1}\right)  .
\]
If $\phi\in\mathcal{D}_{\ast}^{old}\left(  \mathbb{R}\right)  ,$ then
\[
\left\langle \pi^{\prime}\left(  g_{1}\delta_{\ast}\right)  ,\phi\right\rangle
=\phi_{+}\left(  0\right)  .
\]

\end{example}

\subsection{Algebraic and analytic operations in $\mathcal{D}_{\ast,a}%
^{\prime}\left(  \mathbb{R}\right)  \label{subsec3.4}$}

Naturally, we define the algebraic and anlytic operations in $\mathcal{D}%
_{\ast,a}^{\prime}\left(  \mathbb{R}\right)  $ in the same way they are
defined for the usual distributions, namely, by duality.

\subsubsection{Basic definitions}

Let $f,g\in\mathcal{D}_{\ast,a}^{\prime}\left(  \mathbb{R}\right)  ,$
$\phi\left(  x\right)  \in\mathcal{D}_{\ast,a}\left(  \mathbb{R}\right)  ,$
and $\lambda\in\mathbb{C}$. Then $f+\lambda g\in\mathcal{D}_{\ast,a}^{\prime
}\left(  \mathbb{R}\right)  $ is given as
\begin{equation}
\left\langle f+\lambda g,\phi\right\rangle =\left\langle f,\phi\right\rangle
+\lambda\left\langle g,\phi\right\rangle . \label{linearcombination}%
\end{equation}
Let $c\in\mathbb{R}$, translations are handled by
\begin{equation}
\left\langle f\left(  x+c\right)  ,\phi\left(  x\right)  \right\rangle
=\left\langle f\left(  x\right)  ,\phi\left(  x-c\right)  \right\rangle .
\label{translation}%
\end{equation}
Notice that $f\in\mathcal{D}_{\ast,a}^{\prime}\left(  \mathbb{R}\right)  $
while the translation $f\left(  x+c\right)  \in\mathcal{D}_{\ast,a-c}^{\prime
}\left(  \mathbb{R}\right)  .$ Observe that any distribution $g$ of the space
$\mathcal{D}_{\ast,a}^{\prime}\left(  \mathbb{R}\right)  $ can be written as
$g\left(  x\right)  =f\left(  x-a\right)  $ for some $f\in\mathcal{D}_{\ast
}^{\prime}\left(  \mathbb{R}\right)  $, and this justifies studying most
results in $\mathcal{D}_{\ast}^{\prime}\left(  \mathbb{R}\right)  $ only.

Moreover,
\[
\left\langle f\left(  cx\right)  ,\phi\left(  x\right)  \right\rangle
=\frac{1}{\left\vert c\right\vert }\left\langle f\left(  x\right)
,\phi\left(  x/c\right)  \right\rangle .
\]

\bigskip

\subsubsection{Multiplication}

The operation of multiplication is defined by duality: suppose for any
$\phi\in\mathcal{D}_{\ast,a}\left(  \mathbb{R}\right)  ,$ $\psi\phi$ is still
an element in $\mathcal{D}_{\ast,a}\left(  \mathbb{R}\right)  .$ Then if
$\rho\in\mathcal{D}_{\ast,a}^{\prime}\left(  \mathbb{R}\right)  ,$ we define
\[
\left\langle \psi\rho,\phi\right\rangle :=\left\langle \rho,\psi
\phi\right\rangle .
\]
$\psi$ is called a ``multiplier" of $\mathcal{D}_{\ast,a}\left(  \mathbb{R}%
\right)  $ and $\mathcal{D}_{\ast,a}^{\prime}\left(  \mathbb{R}\right)  .$

From the definition one can see that the space of multipliers for a space of
test functions and for its dual space are the same, their Moyal algebra. The
space of multipliers of the spaces of standard test functions and standard
distributions $\mathcal{D}\left(  \mathbb{R}\right)  $ and $\mathcal{D}%
^{\prime}\left(  \mathbb{R}\right)  $ is the space $\mathcal{E}\left(
\mathbb{R}\right)  $ of all smooth function in $\mathbb{R}.$

\begin{definition}
A function $\psi$ defined in $\mathbb{R}\backslash\left\{  a\right\}  $
belongs to $\mathcal{E}_{\ast,a}\left(  \mathbb{R}\right)  $ if $\psi$ is
smooth in $\mathbb{R}\backslash\left\{  a\right\}  $ and if for each ordinary
test function $\rho\in\mathcal{D}\left(  \mathbb{R}\right)  $ the product
$\rho\psi$ belongs to $\mathcal{D}_{\ast,a}\left(  \mathbb{R}\right)  .$
\end{definition}

Clearly, $\mathcal{E}_{\ast,a}\left(  \mathbb{R}\right)  $ contains all
functions that are smooth in $\mathbb{R}\backslash\left\{  a\right\}  $ and
having an asymptotic expansion $\underset{j=m}{\overset{\infty}{\sum}}%
a_{j}\left(  \mathbf{w}\right)  r^{j}$ as $x\rightarrow a.$ Now we show that
$\mathcal{E}_{\ast,a}\left(  \mathbb{R}\right)  $ actually is consisted by
such functions. Indeed, if we pick $\rho\in\mathcal{D}\left(  \mathbb{R}%
\right)  $ such that $\rho\left(  x\right)  \equiv1$ in a neighbourhood of $a,$
then we see that each $\psi\in\mathcal{E}_{\ast,a}\left(  \mathbb{R}\right)  $
actually admits and expansion $\underset{j=m}{\overset{\infty}{\sum}}a_{j}\left(
\mathbf{w}\right)  r^{j}$ as $x\rightarrow a.$

It is clear that we have the following proposition.

\begin{proposition}
The Moyal algebra of $\mathcal{D}_{\ast,a}\left(  \mathbb{R}\right)  $ and of
$\mathcal{D}_{\ast,a}^{\prime}\left(  \mathbb{R}\right)  ,$ namely the space
of multipliers, is $\mathcal{E}_{\ast,a}\left(  \mathbb{R}\right)  .$
\end{proposition}

\begin{example}
The Heaviside function $H\left(  x\right)  =a_{0}\left(  \mathbf{w}\right)
=\left\{
\begin{array}
[c]{c}%
1\text{ \ \ when }\mathbf{w=1}\\
0\text{ \ \ when }\mathbf{w=-1}%
\end{array}
\right.  $ is NOT a multiplier of $D^{\prime}\left(  R\right)  ,$ but it is a
multiplier of $D_{\ast}^{\prime}\left(  R\right)  .$
\end{example}

\subsubsection{Derivatives of thick distributions}

The derivatives of thick distributions are defined also by duality.

\begin{definition}
If $f\in\mathcal{D}_{\ast,a}^{\prime}\left(  \mathbb{R}\right)  $ then its
thick distributional derivative $d^{\ast}f/dx$ is defined as%
\[
\left\langle \frac{d^{\ast}f}{dx},\phi\right\rangle =-\left\langle
f,\frac{d\phi}{dx}\right\rangle ,\ \ \ \phi\in D_{\ast,a}\left(  R\right)  .
\]

\end{definition}

From the discussion in the previous sections, one sees that in general the
spaces $D_{\ast,a}^{\left[  m\right]  }\left(  R\right)  $ are not closed
under differentiation. But $D_{\ast,a}^{\left[  0\right]  }\left(  R\right)  $
is closed under differentiations.

About the notation: $d/dx$ means the ordinary derivative, $\overline{d}/dx$
means the distributional derivative, $d^{\ast}/dx$ means the thick
distributional derivative. Similar to the higher dimensional case, let
$f\in\mathcal{D}_{\ast,a}^{\prime}\left(  \mathbb{R}\right)  $ and
$g=\pi\left(  f\right)  ,$ we have:%

\begin{align*}
\left\langle \pi\left(  \frac{d^{\ast}f}{dx}\right)  ,\phi\right\rangle  &
=\left\langle \frac{d^{\ast}f}{dx},i\left(  \phi\right)  \right\rangle
=-\left\langle f,\frac{di\left(  \phi\right)  }{dx}\right\rangle \\
&  =-\left\langle f,i\left(  \frac{d\phi}{dx}\right)  \right\rangle
=-\left\langle \pi\left(  f\right)  ,\frac{d\phi}{dx}\right\rangle \\
&  =\left\langle \frac{\overline{d}g}{dx},\phi\right\rangle .
\end{align*}
Hence we obtain the following proposition:

\begin{proposition}
Let $f\in\mathcal{D}_{\ast,a}^{\prime}\left(  \mathbb{R}\right)  ,$ then
\[
\pi\left(  \frac{d^{\ast}f}{dx}\right)  =\frac{\overline{d}\pi\left(
f\right)  }{dx}.
\]

\end{proposition}

Because the derivative of test functions satisfies the product rule, the
derivative of thick distrubutions also satisfy the product rule:%
\[
\frac{d^{\ast}\left(  \psi f\right)  }{dx}=\frac{d\psi}{dx}f+\psi\frac
{d^{\ast}f}{dx},\text{ \ \ }f\in\mathcal{D}_{\ast,a}^{\prime}\left(
\mathbb{R}\right)  ,\psi\in\mathcal{E}_{\ast,a}\left(  \mathbb{R}\right)  .
\]
Next let us compute a very important example of thick derivatives, the thick
derivative of the Heaviside function.

\begin{example}
The Heaviside function in $\mathcal{D}_{\ast}^{\prime}\left(  \mathbb{R}%
\right)  $ is defined in the example \ref{HeavisideFunction}, now let us
compute the derivative of it. By definition, for any $\phi\in\mathcal{D}%
_{\ast}\left(  \mathbb{R}\right)  $
\[
\left\langle \frac{d^{\ast}\left(  Pf\left(  H\left(  x\right)  \right)
\right)  }{dx},\phi\right\rangle =-\left\langle Pf\left(  H\left(  x\right)
\right)  ,\frac{d\phi}{dx}\right\rangle .
\]
Now suppose $d\phi/dx$ has asymptotic expansion $d\phi/dx\sim
\underset{j=m}{\overset{+\infty}{\sum}}b_{j}\left(  \mathbf{w}\right)  r^{j}$
and $\phi$ has the asymptotic expansion $\phi\sim
\underset{j=m+1}{\overset{+\infty}{\sum}}a_{j}\left(  \mathbf{w}\right)
r^{j}.$ Then by the definition given in example \ref{HeavisideFunction},%
\begin{align}
-\left\langle Pf\left(  H\left(  x\right)  \right)  ,\frac{d\phi}%
{dx}\right\rangle  &  =-\int_{\left\vert x-a\right\vert \geq A}H\left(
x\right)  \frac{d\phi}{dx}dx+\int_{\left\vert x-a\right\vert <A}H\left(
x\right)  \left(  \frac{d\left(  \phi\left(  x\right)  \right)  }%
{dx}-\underset{j\leq-1}{\sum}b_{j}\left(  w\right)  r^{j}\right)
dx\label{derivative of the Heaviside}\\
&  +\underset{j<-1}{\sum}b_{j}\left(  \mathbf{1}\right)  \frac{A^{j+1}}%
{j+1}+b_{-1}\left(  \mathbf{1}\right)  \log A\nonumber\\
&  =-\int_{A}^{+\infty}\frac{d\phi}{dx}dx-\int_{0}^{A}\left(  \frac{d\left(
\phi\left(  x\right)  \right)  }{dx}-\underset{j\leq-1}{\sum}b_{j}\left(
w\right)  r^{j}\right)  dx\nonumber\\
&  -\underset{j<-1}{\sum}b_{j}\left(  \mathbf{1}\right)  \frac{A^{j+1}}%
{j+1}-b_{-1}\left(  \mathbf{1}\right)  \log A\nonumber\\
&  =\phi\left(  A\right)  -\phi\left(  A\right)  +a_{0}\left(  \mathbf{1}%
\right)  +0=a_{0}\left(  \mathbf{1}\right)  =\left\langle g_{1}\delta_{\ast
},\phi\right\rangle _{\mathcal{D}_{\ast}^{\prime}\left(  \mathbb{R}\right)
\times\mathcal{D}_{\ast}\left(  \mathbb{R}\right)  }\nonumber
\end{align}
Thus the derivative of the Heaviside function is $g_{1}\delta_{\ast},$ which
is defined in the example \ref{Example of thick delta}.

Now consider the projection of the derivative of the Heaviside function
$d^{\ast}\left(  Pf\left(  H\left(  x\right)  \right)  \right)  /dx=g_{1}%
\delta_{\ast}$ onto the usual distribution space $\mathcal{D}^{\prime}\left(
\mathbb{R}\right)  :$%
\begin{equation}
\left\langle \pi\left(  g_{1}\delta_{\ast}\right)  ,\phi\right\rangle
=\left\langle g_{1}\delta_{\ast},i\left(  \phi\right)  \right\rangle
=\phi\left(  0\right)  =\left\langle \delta,\phi\right\rangle .
\label{projg1delta}%
\end{equation}
Keep in mind that $\pi\left(  Pf\left(  H\left(  x\right)  \right)  \right)
=H\left(  x\right)  ,$ the usual Heaviside function. Hence $\pi\left(
d^{\ast}\left(  Pf\left(  H\left(  x\right)  \right)  \right)  /dx\right)
=\delta\left(  x\right)  =\overline{d}\left(  \pi\left(  Pf\left(  H\left(
x\right)  \right)  \right)  \right)  /dx$ as expected.
\end{example}

Next let me present an application of the derivative of thick distributions.

\begin{problem}
Paskusz \cite{Paskusz} pointed out that the following proof is problematic,
where $H\left(  x\right)  $ is the usual Heaviside function:

Since $H\left(  x\right)  =H^{2}\left(  x\right)  $, taking the distributional
derivative on both sides, we have $\delta\left(  x\right)  =2H\left(
x\right)  \delta\left(  x\right)  .$ Hence $H\left(  x\right)  \delta\left(
x\right)  =\frac{1}{2}\delta\left(  x\right)  .$ However, if we multiply
$H\left(  x\right)  $ on both sides again we will get $\frac{1}{2}%
\delta\left(  x\right)  =H\left(  x\right)  \delta\left(  x\right)
=H^{2}\left(  x\right)  \delta\left(  x\right)  =\frac{1}{2}H\left(  x\right)
\delta\left(  x\right)  =\frac{1}{4}\delta\left(  x\right)  ,$ hence we have
$\frac{1}{2}=\frac{1}{4},$ which is clearly wrong.

The key observation of this mistake is that $H\left(  x\right)  \cdot H\left(
x\right)  $ is not a well-defined distribution, that is, $\left\langle
H^{2}\left(  x\right)  ,\phi\left(  x\right)  \right\rangle _{\mathcal{D}%
^{\prime}\left(  \mathbb{R}\right)  \times\mathcal{D}\left(  \mathbb{R}%
\right)  }=\left\langle H\left(  x\right)  ,H\left(  x\right)  \phi\left(
x\right)  \right\rangle $ is not well-defined since $H\left(  x\right)
\phi\left(  x\right)  $ is not a usual test function in $\mathcal{D}\left(
\mathbb{R}\right)  :$ it has a jump discontinuity. Thus we cannot simply apply
the distributional derivative on both sides of the equaion $H\left(  x\right)
=H^{2}\left(  x\right)  .$
\end{problem}

\begin{solution}
In the sense of thick distributions we can restate this whole story in a
rigorous way thus to avoid the mistake $\frac{1}{2}=\frac{1}{4}.$ In fact,
it's easy to see that $H\left(  x\right)  $ is a multiplier of the thick
distributions, i.e. $H\left(  x\right)  f\left(  x\right)  \in\mathcal{D}%
_{\ast}^{\prime}\left(  \mathbb{R}\right)  $ for any $f\left(  x\right)
\in\mathcal{D}_{\ast}^{\prime}\left(  \mathbb{R}\right)  .$ Thus $H\left(
x\right)  \cdot H\left(  x\right)  $ should be viewed as a multiplier times a
thick distribution: $H\left(  x\right)  \cdot Pf\left(  H\left(  x\right)
\right)  $. Then%
\begin{align*}
\frac{d^{\ast}\left(  H\left(  x\right)  \cdot Pf\left(  H\left(  x\right)
\right)  \right)  }{dx}  &  =\frac{d\left(  H\left(  x\right)  \right)  }%
{dx}Pf\left(  H\left(  x\right)  \right)  +H\left(  x\right)  \frac{d^{\ast
}\left(  Pf\left(  H\left(  x\right)  \right)  \right)  }{dx}\\
&  =0+g_{1}\delta_{\ast}%
\end{align*}
from the previous discussions. And equation \ref{projg1delta} tells us that
$\pi\left(  \frac{d^{\ast}\left(  H\left(  x\right)  \cdot Pf\left(  H\left(
x\right)  \right)  \right)  }{dx}\right)  =\delta.$

On the other hand, it is clear that $H\left(  x\right)  \cdot Pf\left(
H\left(  x\right)  \right)  =Pf\left(  H\left(  x\right)  \right)  ,$ taking
derivatives on both sides yields%
\[
g_{1}\delta_{\ast}=g_{1}\delta_{\ast}.
\]

\end{solution}

\nocite{*}
\bibliographystyle{plain}
\bibliography{reconstruct}

\end{document}